\documentclass{amsart}

\usepackage{graphicx}
\usepackage{amsmath, amsfonts, amssymb}

\newcommand{\Z}{\mathbb{Z}}
\newcommand{\R}{\mathbb{R}}

\newcommand{\eps}{\epsilon}
\newcommand{\lagvel}{\nu}

\DeclareMathOperator{\im}{Im}
\DeclareMathOperator{\re}{Re}
\DeclareMathOperator{\sgn}{sgn}

\renewcommand\L{\mathcal{L}}
\newcommand\What{\widehat{W}}
\newcommand\Shat{\widehat{S}}
\newcommand\uhat{\widehat{u}}
\newcommand\vhat{\widehat{v}}

\begin{document}

\title[Reflection of sound waves off an entropy wave]{On the resonant reflection of weak, nonlinear sound waves off an entropy wave}

\author{John K. Hunter}
\address{Department of Mathematics, University of California at Davis}
\email{jkhunter@ucdavis.edu}
\thanks{JKH was supported by the NSF under grant number DMS-1616988}
\author{Evan B. Smothers}
\address{Uber Technologies, Inc.,  San Francisco}
\email{sevan@uber.com}

\date{October 14, 2018}

\begin{abstract}
We derive a degenerate quasilinear Schr\"odinger equation that describes the resonant reflection of very weak, nonlinear sound waves
off a weak sawtooth entropy wave.
\end{abstract}

\maketitle

\section{Introduction}

Away from a vacuum, the Eulerian form of the one-dimensional, nonisentropic gas dynamics equations is a $3\times 3$ strictly hyperbolic system of conservation laws \cite{whitham}. There are three wave fields: genuinely nonlinear left and right moving sound waves, and linearly degenerate entropy waves.
In the absence of an entropy wave, the interaction between weakly nonlinear, counter-propagating sound waves is negligible to leading order in the wave amplitudes. Asymptotic expansions show that the velocity perturbations in the sound waves satisfy separate inviscid Burgers equations, and, after a sufficiently long time, the left and right moving sound waves form shocks and decay in the same way as they would in a unidirectional wave.

The presence of an entropy wave can have a dramatic effect on the dynamics of spatially periodic sound waves. When the wavenumber of the entropy wave is twice that of the sound waves, the sound waves reflect resonantly off the entropy wave, leading to temporal oscillations in the sound waves. The resulting alternation between compression and expansion can delay
the formation of shocks and weaken them greatly due to the absorption of expansion waves
created by the oscillations.

The interaction of weak sound waves and entropy waves with amplitudes of the same order $\epsilon \ll 1$ was studied by Majda, Rosales and Schonbeck \cite{mrs}.
In that case, the normalized amplitudes $u(x,t)$, $v(x,t)$ of the left and right moving sound waves, in separate reference frames moving with their linearized sound speeds, satisfy a system of the form
\begin{align}
\begin{split}
&u_t + \left(\frac{1}{2}u^2\right)_x = \frac{1}{2\pi}\int_{-\pi}^{\pi} K(x-y) v_y(y,t)\, dy,
\\
&v_t + \left(\frac{1}{2}v^2\right)_x = \frac{1}{2\pi}\int_{-\pi}^{\pi} K(y-x) u_y(y,t)\, dy.
\end{split}
\label{gmrs}
\end{align}
Here, the kernel $K(x)$ is a rescaled entropy-wave profile, which is not affected by the sound waves to leading order,
and all functions are $2\pi$-periodic in $x$ with zero mean.
The convolution terms in \eqref{gmrs} describe the resonant reflection of the left and right moving sound waves off the entropy wave, while the Burgers terms describe the nonlinear deformation of the sound waves.
The validity of these asymptotic equations for weak solutions is proved in \cite{schochet}, and the global existence of entropy solutions of \eqref{gmrs} is proved in \cite{xin}.

In general, the system \eqref{gmrs} is spatially nonlocal and  dispersive. One of the simplest cases, and one with particularly interesting dynamics, is when  the entropy wave is a sawtooth wave
\[
K(x) = \begin{cases} x+\pi &\text{for $-\pi<x<0$,}
\\
 x-\pi &\text{for $0<x<\pi$,}
 \end{cases}
\qquad K(x+2\pi) = K(x),
\]
which leads to the local, nondispersive system \cite{mrs}
\begin{align}
\begin{split}
&u_t +\left(\frac{1}{2}u^2\right)_x = v, \qquad
v_t + \left(\frac{1}{2}v^2\right)_x = -u.
\end{split}
\label{MRS}
\end{align}
In this case, the linearized sound waves oscillate with frequency $1$
independent of their wavenumber and have zero group velocity.
We refer to \eqref{MRS} as the Majda-Rosales-Schonbeck (MRS) equation.

In this paper, we consider the resonant reflection of sound waves off an entropy wave when the sound waves are much weaker than the entropy wave. Specifically,  we suppose that the entropy wave is a sawtooth wave with amplitude of the order $\epsilon\ll 1$ and the sound waves have amplitudes of the order
$\epsilon^{3/2} \ll \epsilon$. This particular choice of scaling leads to a dominant balance between weak nonlinearity and weak dispersion.

There are then three timescales:
on times of the order $1$, the sound waves propagate at the linearized sound speed; on times of the order
$1/\epsilon$, the sound wave profiles oscillate linearly with constant frequency as a result of their reflection off the entropy wave; and on times of the order $1/\epsilon^2$, the spatial profiles of these oscillations are modulated by the nonlinear self-interaction of the sound waves.

It is convenient to carry out the asymptotic expansion in Lagrangian coodinates, when the entropy wave is independent of time. As summarized in Section~\ref{sec:gd}, the deformation $x\mapsto \chi(x,t)$ of the compressible fluid satisfies the quasilinear wave equation \eqref{lagd} with wave speed \eqref{c0exp}.

 In suitably nondimensionalized variables, our asymptotic solution for the Lagrangian velocity $\lagvel = \chi_t$ of the fluid is given by
\begin{equation}
\lagvel(x,t;\epsilon) = \epsilon^{3/2}\left\{ a(t-x,\epsilon^2 t) -ia(t + x,\epsilon^2t) \right\}e^{-i\epsilon t}
 + \text{c.c.} + O(\epsilon^2)
 \label{lagvel}
\end{equation}
as $\epsilon \to 0^+$, where the complex-valued amplitude function $a(z,\tau)$ is $2\pi$-periodic in $z$ with zero mean and satisfies the degenerate quasilinear Schr\"odinger equation
\begin{align}
&i\left(a_{\tau} + \frac{1}{2}\partial_z^{-1}a - \frac{\pi^2}{2} a_z \right)= \frac{(\gamma+1)^2}{6} (|a|^2 a_z)_z.
\label{asyeq}
\end{align}
Equation \eqref{lagvel} describes right and left moving sound waves with linearized speed $1$ and complex amplitude functions $a$ and $-ia$, respectively. The sound waves oscillate slowly with frequency $\epsilon$ and evolve over longer times of the order $1/\epsilon^{2}$ according to \eqref{asyeq}.
As described in Section~\ref{sec:mrs}, the dispersionless version of \eqref{asyeq}, given in \eqref{dqnls},
also arises directly from the MRS equation \eqref{MRS}.

When $|a|$ is bounded away from zero, one expects that the qualitative behavior of \eqref{asyeq} is similar to a Schr\"odinger equation. However, energy estimates for the derivatives of $a$ fail when $a\to 0$.
Numerical solutions of \eqref{asyeq} shown in Section~\ref{sec:num_gd} suggest that $a$ may become zero at some point, even if it is initially nonzero, leading to the generation of smaller-scale waves and a possible loss of smoothness in the solution.

From a more general perspective, the much richer dynamical possibilities for spatially periodic solutions of the nonisentropic gas dynamics equations in comparison with the isentropic equations is the simplest and most physically important example of the different behavior of large-variation solutions of $n\times n$ nonlinear hyperbolic conservation laws for $n \ge 3$ in comparison with $n=2$.

Glimm \cite{glimm} proved the existence of global weak entropy solutions of $n\times n$ systems of strictly hyperbolic conservation laws in one space dimension for $u : \R\times \R^+\to \R^n$,
\[
u_t + f(u)_x = 0\qquad f : \R^n\to \R^n,
\]
with genuinely nonlinear or linearly degenerate wave fields and initial data with sufficiently small total variation. These Glimm solutions were later shown to be unique by Bianchini and Bressan \cite{bressan}. In this case, the different wave fields can undergo only a limited amount of significant interactions, and they eventually separate
into nondecaying linearly degenerate waves and decaying genuinely nonlinear $N$-waves \cite{liu1, liu2}.

For $2\times 2$ genuinely nonlinear systems, Glimm and Lax \cite{glimmlax} proved that spatially periodic solutions that are small in $L^\infty$ (but have infinite total variation on $\R$) decay to zero as $t\to\infty$.
Thus, as a consequence of shock-formation and the resulting decay induced by shocks, the long-time dynamics of solutions is trivial in both of these cases

For $n\ge 3$, the behavior of solutions of $n\times n$ hyperbolic systems with small $L^\infty$-norm and large total variation may be very different, since significant wave interactions can persist for all $t >0$. In particular, resonant three-wave interactions may weaken shocks, or delay their formation, or conceivably prevent their formation entirely \cite{temple}. This may lead to solutions that do not decay to zero as $t\to \infty$ and have nontrivial long-time nonlinear dynamics.
The formal asymptotic solutions for resonant gas dynamics support these suggestions \cite{pego, shefter}, as do the ones derived here.
However, it should be noted that these asymptotic solutions
only apply for long times of some order $1/\epsilon^{n}$
and need not remain valid as $t\to\infty$, leading to difficult KAM-type issues \cite{temple1}.

An outline of this paper is as follows. In Section~\ref{sec:mrs}, we derive the asymptotic equation \eqref{dqnls} for the MRS equation. In Section~\ref{sec:gd}, we summarize the Lagrangian gas dynamics equations, and in Section~\ref{sec:der}, we derive the asymptotic equation \eqref{asyeq}. Finally, in Section~\ref{sec:num}, we present some numerical solutions.

\section{The MRS Equation}\label{sec:mrs}

As a model problem, we first consider the MRS equation \eqref{MRS}, which provides an approximate description of the
resonant reflection of sound waves off a sawtooth entropy wave.
We will show that this system has an asymptotic solution of the form
\begin{align}
\begin{split}
u(x,t;\eps) &= \eps a(x,\eps^2 t) e^{-it} + \text{c.c.} + O(\eps^2),
\\
v(x,t;\eps) &= -i\eps a(x,\eps^2 t) e^{-it}  + \text{c.c.}   + O(\eps^2),
\end{split}
\label{mrsasysol}
\end{align}
as $\epsilon \to 0$, valid on timescales $t=O(\epsilon^{-2})$, where the complex-valued amplitude $a(x,\tau)$ satisfies
a degenerate, quasilinear Schr\"odinger equation \cite{hrev}
\begin{equation}
i a_\tau +\left( \frac{2}{3} |a|^2 a_x\right)_x = 0.
\label{dqnls}
\end{equation}

Equation \eqref{mrsasysol} describes a solution with an arbitrary spatial profile and amplitude of the order $\epsilon$ that oscillates in time with frequency $1$.  Nonlinear effects lead to a slow deformation of the profile over long times of the order $1/\epsilon^2$, according to \eqref{dqnls}.

Using the method of multiple scales, we look for asymptotic solutions  of \eqref{MRS}  of the form
\begin{align*}
u(x,t;\eps) &= \eps u_1(x,t,\eps^2 t) + \eps^2 u_2(x,t,\eps^2 t) + \eps^3 u_3(x,t,\eps^2 t) + O(\epsilon^4),
\\
v(x,t;\eps) &= \eps v_1(x,t,\eps^2 t) + \eps^2 v_2(x,t,\eps^2 t) + \eps^3 v_3(x,t,\eps^2 t) + O(\epsilon^4).
\end{align*}
Substituting these expansions into \eqref{MRS}, expanding derivatives, and equating coefficients of $\epsilon$, $\epsilon^2$, and $\epsilon^3$, we find that the functions $u_j(x,t,\tau)$, $v_j(x,t,\tau)$ with $j=1,2,3$ satisfy
\begin{align}
&u_{1t} = v_1,\qquad v_{1t} = -u_1
\label{pertmrs1}
\\
&u_{2t} + \left(\frac{1}{2} u_1^2\right)_x = v_2,\qquad v_{2t} + \left(\frac{1}{2} v_1^2\right)_x= -u_2,
\label{pertmrs2}
\\
&u_{3t} + u_{1\tau} +  \left(u_1 u_2\right)_x = v_3,\qquad v_{3t} + v_{1\tau} + \left(v_1 v_2\right)_x= -u_3.
\label{pertmrs3}
\end{align}

The solution of \eqref{pertmrs1} for $(u_1,v_1)$ is
\begin{equation}
u_1(x,t,\tau) = a(x,\tau) e^{-it} + \mathrm{c.c.},\quad v_1(x,t,\tau) = -ia(x,\tau) e^{-it} + \mathrm{c.c.},
\label{sol1}
\end{equation}
where $a(x,\tau)$ is an arbitrary complex-valued function.
The solution of \eqref{pertmrs2}  for $(u_2,v_2)$ is then
\begin{align}
\begin{split}
&u_2(x,t,\tau) = B(x,\tau) e^{-2it} - M(x,\tau) + \mathrm{c.c.},
\\
&v_2(x,t,\tau) = C(x,\tau) e^{-2it} + M(x,\tau)+\mathrm{c.c.},
\\
&B = -\left(\frac{1+2i}{6} \right) \left(a^2\right)_x,\quad
C = -\left(\frac{1-2i}{6} \right) \left(a^2\right)_x,
\\
&M = \frac{1}{2}\left(|a|^2\right)_x.
\end{split}
\label{sol2}
\end{align}
where we omit a homogeneous solution that does not enter into the final result.

The solvability condition for the removal of secular terms from solutions $(u_3,v_3)$
of the system
\[
u_{3t} - v_3 = Fe^{-it},\qquad v_{3t} + u_3 = Ge^{-it}
\]
is $F + i G = 0$. We use \eqref{sol1}--\eqref{sol2} in \eqref{pertmrs3}, compute the coefficients of $e^{-it}$, and impose this solvability condition, which gives equation \eqref{dqnls} for $a(x,\tau)$.

We remark that one can derive similar asymptotic solutions for more general balance laws for $u(x,t)\in \R^n$ of the form
\[
u_t + \sum_{j=1}^d f^j(u)_{x^j} = g(u),
\]
where $\nabla_u f^j(0) = 0$, so that the fluxes $f^j : \R^n\to \R^n$ are at least quadratically nonlinear, and
$\nabla_u g(0)$ has purely imaginary eigenvalues, so that the linearized equation has oscillatory solutions. See \cite{smothers}
for more details in the one-dimensional case $d=1$.

The MRS system \eqref{MRS} can be compared with the Burgers-Hilbert equation
\[
u_t + \left(\frac{1}{2} u^2\right)_x = H[u],
\]
where $H$ is the Hilbert transform with symbol $-i\sgn k$.
This equation arises as a description of waves on a vorticity discontinuity and
models unidirectional, real wave motions whose linearized frequency is independent of their wavenumber \cite{biello, bh}.

Since $H^2 = -I$, the functions
$(u,v)$ with $v=H[u]$ satisfy
\[
u_t + \left(\frac{1}{2} u^2\right)_x = v,\qquad
v_t +|\partial_x|\left(\frac{1}{2} u^2\right) = -u,
\]
where $|\partial_x| = H \partial_x$ has symbol $|k|$.
This system has the same linear terms as the MRS equation, but has a nonlocal nonlinear term.
Instead of \eqref{dqnls}, the complex amplitude $a(x,\tau)$ of weakly nonlinear solutions
of the Burgers-Hilbert equation satisfies a nonlocal, degenerate, cubically quasilinear Schr\"odinger-type equation
\[
i a_\tau + P\left[|a|^2 a_x - i a|\partial_x| \left(|a|^2 \right)\right]_x = 0,
\]
where $P$ denotes the projection onto positive wavenumber components, and $Pa=a$ \cite{biello}.

The degenerate quasilinear Schr\"odinger (DQS) equation \eqref{dqnls} is of interest in its own right as a degenerate dispersive analog of the degenerate diffusive porous medium equation $u_t = (u^2u_x)_x$ \cite{pm}. We will study the DQS equation in more detail elsewhere \cite{hs}, but,  by analogy with wetting fronts for the porous medium equation, it is natural to ask about the propagation of finite-speed dispersive fronts into the zero solution for the DQS equation. In Section~\ref{sec:num_mrs}, we show a numerical comparison of front propagation for the MRS and DQS equations. Front propagation in a different degenerate quasilinear Schr\"odinger equation and related degenerate quasilinear KdV equations is analyzed in \cite{germain, germain1}.

\section{Gas Dynamics}\label{sec:gd}
We consider the one-dimensional flow of an inviscid, ideal gas with equation of state
\[
p = K \rho^\gamma e^{s/C_V},
\]
where
$p$, $\rho$, $s$ are the pressure, spatial density, and entropy, respectively, and $K$, $C_V$, $\gamma$ are
positive constants.

The Lagrangian form of the
compressible Euler equations for the deformation $x\mapsto \chi(x,t)$ of a Lagrangian particle $x$
can be written as the quasilinear wave equation (see e.g., \cite{shkoller})
\begin{equation}
\chi_{tt} + \frac{1}{\gamma} c_0^2(x) \left[\frac{1}{\chi_x^{\gamma}}\right]_x = 0,
\label{lagd}
\end{equation}
where we use a reference configuration in which $p=p_0$ is constant
when $\chi(x,t) =x$, and the sound speed $c_0(x)$ is given in terms of the Lagrangian density $\rho_0(x)$ by
\[
c_0(x) = \left(\frac{\gamma p_0}{\rho_0(x)}\right)^{1/2}.
\]

We use nondimensionalized $(x,t)$-variables in which the mean sound speed is equal to $1$
and the perturbations in the sound speed have wavelength
$\pi$ with respect to the Lagrangian variable $x$.
Introducing a small parameter $\epsilon \ll 1$, we suppose that
\begin{equation}
c_0^2(x;\epsilon) = 1+ \epsilon S(2x),
\label{c0exp}
\end{equation}
where the $2\pi$-periodic sawtooth function $S(x)$ is given by
\begin{equation}
S(x) = \begin{cases}
       2(x+\pi) & \text{if $-\pi \le x< 0$,} \\
       2(x-\pi) &\text{if $0<x\le\pi$,}
     \end{cases}\qquad S(x+2\pi) = S(x).
\label{defS}
\end{equation}
The corresponding entropy perturbation is
\begin{equation}
s(x) = s_0 + \gamma C_V \log\left[1 + \epsilon S(2x)\right] = s_0 + \epsilon\gamma C_V S(2x) + O(\epsilon^2).
\label{entexp}
\end{equation}
One could include $O(\epsilon^2)$ perturbations in the expansion \eqref{c0exp} of the wave speed, but for simplicity we assume they are zero. The effect of such perturbations is discussed further at the end of Section~\ref{MRSlin}.

We consider weakly nonlinear asymptotic solutions of \eqref{lagd} which are
$2\pi$-periodic in $x$. In that case, the phases of the left and right moving sound waves and the entropy wave satisfy
$(t+x) - (t-x) = 2x$,
so the sound waves reflect resonantly into each other off the entropy wave.
In the next section, we derive the asymptotic solution of \eqref{lagd}
for the Lagrangian velocity $\lagvel = \chi_t$ that is given in \eqref{lagvel}--\eqref{asyeq}.

\section{Derivation of the asymptotic equation}\label{sec:der}

Writing $\chi(x,t) = x + w(x,t)$ in \eqref{lagd}, we get that the displacement $w$ satisfies
\begin{equation}
w_{tt} + \frac{1}{\gamma} c_0^2(x) \left[\frac{1}{(1 + w_x)^{\gamma}}\right]_x = 0,
 \label{gasdynamicsPDE}
\end{equation}
where $c_0^2$ is given by \eqref{c0exp}.
We look for asymptotic solutions of  \eqref{gasdynamicsPDE} of the form
\begin{align}
\begin{split}
w(x,t;\epsilon) &=  \epsilon^{3/2} W(x,t,T,\tau;\epsilon),\quad T = \epsilon t, \quad\tau = \epsilon^2 t,
\\
W &= W_0 + \epsilon^{1/2} W_1 + \epsilon W_2 + \epsilon^{3/2} W_3 + \epsilon^2 W_4 + O(\epsilon^{5/2}),
\end{split}
\label{ansatz}
\end{align}
where $W(x,t,T,\tau;\epsilon)$ is a $2\pi$-periodic function of $x$ with zero mean.

Using \eqref{c0exp} and \eqref{ansatz} in \eqref{gasdynamicsPDE}, expanding time derivatives, Taylor expanding the
result,  and equating coefficients of powers of $\epsilon^{1/2}$,
we obtain the perturbation equations
\begin{align}
\L W_0 &= 0,\label{pert0} \\
\L W_1 &= 0,\label{pert1}\\
\L W_2 &= S(2x) W_{0xx} - 2 W_{0tT},\label{pert2} \\
\L W_3 &= S(2x) W_{1xx} - 2 W_{1tT} - (\gamma+1)W_{0xx}W_{0x},\label{pert3}  \\
\L W_4 &= S(2x) W_{2xx} - 2W_{2tT} - (\gamma+1)(W_{0x}W_{1x})_x\nonumber
\\
&\qquad\qquad - W_{0TT} - 2W_{0t\tau}\label{pert4},
\end{align}
where $\L = \partial_t^2 - \partial_x^2$.

\subsection{The $W_0$, $W_1$ equations}
We make a change of independent variables $(x,t) \mapsto (\xi,\eta)$, where
\[
\xi = t-x, \qquad \eta = t+x,\qquad \L= 4\partial_{\xi} \partial_{\eta}.
\]
Then the solutions of \eqref{pert0}--\eqref{pert1} are
\begin{align}
\begin{split}
W_0 &= U_0(\xi,T,\tau) + V_0(\eta,T,\tau),
\\
W_1 &= U_1(\xi,T,\tau) + V_1(\eta,T,\tau),
\end{split}
\label{W0W1}
\end{align}
where $U_0$, $U_1$ or $V_0$, $V_1$ are arbitrary $2\pi$-periodic, zero-mean functions of $\xi$ or $\eta$, respectively.

For $j=0,1,2,\dots$, we use the notation
\[
u_j(\xi,T,\tau) = U_{j\xi}(\xi,T,\tau), \qquad v_j(\eta,T,\tau) = V_{j\eta}(\eta,T,\tau),
\]
so that
\begin{equation}
U_{j t} = u_j,\qquad U_{jx} = - u_j,\qquad V_{jt} = v_j,\qquad V_{jx} = v_j.
\label{UVxt}
\end{equation}
The Lagrangian velocity $\lagvel = w_t$ is given to leading order by
\begin{equation}
\lagvel = \epsilon^{3/2}\left(u_0 + v_0\right) + O(\epsilon^2).
\label{lagvelpert}
\end{equation}

\subsection{The $W_2$ equation}
Since $\eta-\xi = 2x$, equation \eqref{pert2} becomes
\begin{align}
4 W_{2\xi\eta} &= S(\eta-\xi)[u_{0\xi} + v_{0\eta}] - 2[u_{0T} + v_{0T}] .
\label{W2equation}
\end{align}
In order for \eqref{W2equation} to have a solution for $W_2$ that is $2\pi$-periodic in $\xi$, $\eta$, the right-hand side must have zero mean with respect to $\xi$, $\eta$. We use $\langle f \rangle^{\xi}$ to denote the average of a $2\pi$-periodic function $f(\xi,\eta)$ with respect to $\eta$ and $\langle f \rangle^{\eta}$ to denote the average with respect to $\xi$,
\[
\langle f \rangle^{\xi} = \frac{1}{2\pi} \int_{-\pi}^{\pi} f(\xi,\eta)\, d\eta,\qquad
\langle f \rangle^{\eta} = \frac{1}{2\pi} \int_{-\pi}^{\pi} f(\xi,\eta)\, d\xi.
\]

Using the fact that $S'(x) = 2-4\pi\delta(x)$, together with
the zero-mean periodicity of $S$, $u_0$, $v_0$,
we compute that
\begin{align}
\begin{split}
&\langle S(\eta-\xi) u_{0\xi}(\xi) \rangle^{\eta} = -2 u_0(\eta),
\qquad
\langle S(\eta-\xi) v_{0\eta}(\eta) \rangle^{\eta} = 0,
\\
&\langle S(\eta-\xi) u_{0\xi}(\xi) \rangle^{\xi}
= 0,
\qquad\qquad
\langle S(\eta-\xi) v_{0\eta}(\eta) \rangle^{\xi}
= 2 v_0(\xi),
\end{split}
\label{Sav}
\end{align}
where we suppress the dependence on $T$ and $\tau$.

Taking the averages of \eqref{W2equation} with respect to $\xi$ and $\eta$, we therefore obtain that
\begin{align}
&u_{0T} = v_0,\qquad
v_{0T} = -u_0.
\label{mrspert}
\end{align}
The solution of \eqref{mrspert} is
\begin{align}
\begin{split}
u_0(\xi,T,\tau) &= a(\xi,\tau) e^{-iT} + \text{c.c.}, \\
v_0(\eta,T,\tau) &= -ia(\eta,\tau) e^{-iT} + \text{c.c.},
\end{split}
\label{mrspertsol}
\end{align}
where $a(z,\tau)$ is an arbitrary complex-valued function that is $2\pi$-periodic in $z$ with zero mean.
The use of \eqref{mrspertsol} in \eqref{lagvelpert} gives \eqref{lagvel}.

Using \eqref{Sav}--\eqref{mrspert}, we can rewrite equation \eqref{W2equation} for $W_2$ as
\begin{equation}
4 W_{2\xi\eta} = S(\eta-\xi)u_{0\xi} - \langle S(\eta-\xi) u_{0\xi} \rangle^{\eta} + S(\eta-\xi)v_{0\eta} - \langle S(\eta-\xi) v_{0\eta}\rangle^{\xi}. \label{W2withaverages}
\end{equation}
We Fourier expand the functions in this equation as
\begin{align}
\begin{split}
&u_{0\xi}(\xi) = \sum_{k \in \mathbb{Z}} \uhat_k \ e ^{i k \xi}, \qquad\qquad\qquad\quad v_{0 \eta}(\eta) = \sum_{l \in \mathbb{Z}} \vhat_l e^{i l \eta}, \\
&W(\xi,\eta) = \sum_{m,n \in \mathbb{Z}} \What_{mn} \ e^{i(m\xi + n \eta)}, \qquad S(\eta-\xi) = \sum_{p \in \mathbb{Z}} \Shat_p e^{i p (\eta-\xi)}.
\end{split}
\label{fourierexp}
\end{align}
Since $\Shat_0 = 0$, it follows that
\begin{align*}
S(\eta-\xi)u_{0\xi} - \langle S(\eta-\xi) u_{0\xi} \rangle^{\eta}
&= \sum_{m,n \in \Z_*} \Shat_{n}\uhat_{m+n} e^{i(m\xi+n\eta)},
\\
S(\eta-\xi)v_{0\eta} - \langle S(\eta-\xi) v_{0\eta}\rangle^{\xi}
&=\sum_{m,n \in \Z_*} \Shat_{-m}\vhat_{m+n} e^{i(m\xi+n\eta)},
\end{align*}
where $\Z_* = \Z\setminus\{0\}$.
Hence, the Fourier coefficients of a particular solution of the nonhomogeneous equation \eqref{W2withaverages} are given by
\begin{align*}
\What_{mn} &= -\frac{\Shat_{n}\uhat_{m+n}+\Shat_{-m}\vhat_{m+n}}{4mn}
\end{align*}
for $m,n \neq 0$, and the general solution of \eqref{W2withaverages}  is
\begin{align}
\begin{split}
W_2(\xi,\eta) &= U_2(\xi) + V_2(\eta)
\\
&\quad - \sum_{m,n \in \Z_*} \frac{\Shat_{n}\uhat_{m+n}+\Shat_{-m}\vhat_{m+n}}{4mn} e^{i(m\xi+n\eta)},
\end{split}
\label{W2sol}
\end{align}
where $U_2$, $V_2$ are arbitrary functions of integration.

\subsection{The $W_3$ equation}
Next, we examine equation \eqref{pert4} for $W_3$, which can be written as
\begin{align}
\begin{split}
4 W_{3\xi\eta} &=  -(\gamma+1)[u_{0\xi} + v_{0\eta}][v_0-u_0]
\\
&\quad + S(\eta-\xi)[u_{1\xi} + v_{1\eta}] - 2[u_{1T} + v_{1T}].
\end{split}
\label{w3eq}
\end{align}
As with $W_2$, we impose the solvability conditions that follow by averaging the equation with respect to $\xi$ and $\eta$.
Using \eqref{Sav} and the fact that $\langle F G_{\eta} \rangle^{\xi} = 0$, $\langle F_{\xi} G \rangle^{\eta} = 0$ for any functions $F(\xi)$, $G(\eta)$, we find that $u_1$ and $v_1$ satisfy
\begin{align*}
u_{1T} - v_1 &=  \frac{\gamma+1}{2} u_{0\xi}u_0, \qquad
v_{1T} + u_1 = - \frac{\gamma+1}{2} v_{0\eta}v_0.
\end{align*}
Using \eqref{mrspertsol} in these equations and computing the solution, we get that
\begin{align}
\begin{split}
u_1 &= \frac{\gamma+1}{6}(-1+2i) aa_x e^{-2iT} - \frac{\gamma+1}{4} |a|^2_z + \text{c.c.}\\
v_1 &= \frac{\gamma+1}{6}(1+2i) a a_x e^{-2iT} - \frac{\gamma+1}{4} |a|^2_z + \text{c.c.},
\end{split}
\label{mrspertsol1}
\end{align}
where we omit a solution of the homogeneous equation, which does not enter into the final result. We can then solve
\eqref{w3eq} for $W_3$, but we do not need an explicit expression since $W_3$  does not appear in the equation for $W_4$.

\subsection{The $W_4$ equation}
We derive an equation for $a(z,\tau)$ by imposing the previous solvability conditions on equation \eqref{pert4} for $W_4$, which has the form
\begin{align*}
4 W_{4\xi\eta} &= f
\\
f&=-(\gamma+1)(W_{0x}W_{1x})_x + S(2x) W_{2xx} - 2 W_{2tT} - W_{0TT} - 2 W_{0t\tau}.
\end{align*}

The first solvability condition is
\begin{equation}
\langle f \rangle^{\xi} = 0,\qquad  \langle f \rangle^{\eta} = 0,
\label{fsolcon}
\end{equation}
which ensures that there is a solution for $W_4$ that is $2\pi$-periodic in $(\xi,\eta)$. We
compute the averages of each of the terms in $f$ separately.

Using \eqref{W0W1}--\eqref{UVxt}, we have
\begin{align*}
&\langle W_{0TT}  \rangle^{\xi} = U_{0TT},
\qquad
\langle V_{0TT} \rangle^{\eta} = V_{0TT},
\\
&\langle W_{0t\tau} \rangle^{\xi} = u_{0\tau},
\qquad\quad
\langle W_{0t\tau} \rangle^{\eta} = v_{0\tau}.
\end{align*}
In addition,
\begin{align*}
(W_{0x}W_{1x})_x &= (\partial_{\eta} - \partial_{\xi})[(v_0-u_0)(v_1-u_1)],
\end{align*}
so that
\begin{align*}
&\langle (W_{0x}W_{1x})_x \rangle^{\xi} = - (u_0u_1)_{\xi},\qquad
\langle (W_{0x}W_{1x})_x \rangle^{\eta} = (v_0v_1)_{\eta}.
\end{align*}
From \eqref{W2sol} and \eqref{UVxt} we also obtain that
\begin{align*}
\langle W_{2tT} \rangle^{\xi} &= u_{2T},\qquad
\langle W_{2tT} \rangle^{\eta} = v_{2T}.
\end{align*}

To compute the averages of $S(2x) W_{2xx}$, we use the previous Fourier expansions, which gives
\begin{align}
\begin{split}
&S(2x) W_{2xx} = S(\eta-\xi)[u_{2\xi} + v_{2\eta}] \\
&+\left(\sum_{p \in \Z} \Shat_p e^{i p (\eta-\xi)} \right)\left(\sum_{m,n\in\Z_*} \frac{(n-m)^2}{4mn}[\Shat_{n}\uhat_{m+n}+\Shat_{-m}\vhat_{m+n}]e^{i(m\xi+n\eta)}\right) \label{Sumformula}
\end{split}
\end{align}
In the $\eta$ average, we retain only the $\eta = 0$ Fourier coefficients, which corresponds to taking $p+n = 0$, and in the $\xi$ average, we retain only the $\xi = 0$ Fourier coefficients, which corresponds to taking $m-p = 0$.  Applying this observation to \eqref{Sumformula}, rewriting the result as sum over $k=m+n$, and using \eqref{Sav} for the first term, we get that
\begin{align}
&\langle S(\eta-\xi) [\partial_{\eta}-\partial_{\xi}]^2 W_2 \rangle^{\xi} = 2v_2(\xi)
\nonumber\\
&\qquad+ \sum_{k \in\Z}\left(\sum_{n\in\Z\setminus\{0,k\}}\frac{(2n-k)^2}{4(k-n)n}\Shat_{-n}[\Shat_n\uhat_k + \Shat_{n-k} \vhat_k]\right)e^{ik\xi}
\label{Sxiav}\\
&\langle S(\eta-\xi) [\partial_{\eta}-\partial_{\xi}]^2 W_2 \rangle^{\eta} = - 2 u_2(\eta)
\nonumber\\
&\qquad+ \sum_{k \in\Z}\left(\sum_{n\in\Z\setminus\{0,k\}}\frac{(2n-k)^2}{4(k-n)n}\Shat_{k-n}[\Shat_n\uhat_k + \Shat_{n-k} \vhat_k]\right)e^{ik\eta}.
\label{Setaav}
\end{align}

For the sawtooth wave \eqref{defS}, we have
\begin{align*}
\Shat_k &= \frac{1}{2\pi} \int_{-\pi}^{\pi} S(x) e^{-i k x } \ dx = \frac{2i}{k}.
\end{align*}
Using this value in \eqref{Sxiav}, computing the resulting sums, and using  \eqref{fourierexp}, we find that
\begin{align*}
\langle S(\eta-\xi) [\partial_{\eta}-\partial_{\xi}]^2 W_2 \rangle^{\xi} &= 2v_2(\xi) - \sum_{k \in \Z_*}\left\{ \pi^2 \uhat_k + \left(\frac{2\pi^2}{3} + \frac{2}{k^2}\right) \vhat_k \right\}e^{i k \xi} \\
&= 2v_2(\xi) - \pi^2 u_{0\xi}(\xi) - \frac{2\pi^2}{3} v_{0\xi}(\xi) - 2\partial_\xi^{-1} v_0(\xi),
\end{align*}
where $\partial_\xi^{-1}$ is the Fourier multiplier operator on periodic, zero-mean functions with symbol ${1}/{ik}$.
A similar computation in \eqref{Setaav} gives
\begin{align*}
&\langle S(\eta-\xi) [\partial_{\eta}-\partial_{\xi}]^2 W_2 \rangle^{\eta}
\\
&\qquad = -2u_2(\eta) - \pi^2 v_{0\eta}(\eta) - \frac{2\pi^2}{3} u_{0\eta}(\eta) - 2\partial_\eta^{-1}u_0(\eta).
\end{align*}

Combining these results, we find that the solvability condition \eqref{fsolcon} is
\begin{align}
\begin{split}
u_{2T} - v_2 &= \frac{\gamma+1}{2} (u_0 u_1)_{\xi} - \frac{\pi^2}{2} u_{0\xi} - \frac{\pi^2}{3} v_{0\xi} - \partial_\xi^{-1} v_0 - \frac{1}{2}  \partial_\xi^{-1}u_{0TT} - u_{0\tau}, \\
v_{2T} + u_2 &= -\frac{\gamma+1}{2} (v_0 v_1)_{\eta} - \frac{\pi^2}{2} v_{0\eta} - \frac{\pi^2}{3} u_{0\eta} - \partial_\eta^{-1}u_0 - \frac{1}{2} \partial_\eta^{-1}v_{0TT} - v_{0\tau}.
\end{split}
\label{u2v2eq}
\end{align}

The first equation is an equation in $\xi$ and the second is an equation in $\eta$.  We replace $\xi$ by $z$ in the first equation and $\eta$ by $z$ in the second equation to get a system  for functions of $z$.

We now impose the second solvability condition, which comes from the requirement that \eqref{u2v2eq} has solutions for $(u_2,v_2)$ that are $2\pi$-periodic in $T$, to avoid secular terms in $T$. If
\begin{align*}
&u_{2T} - v_2 = Fe^{-iT} + \text{c.c.} + \text{n.r.t.},
\qquad
v_{2T} + u_2 = Ge^{-iT} + \text{c.c.} + \text{n.r.t.},
\end{align*}
where $\text{n.r.t.}$ stands for non-resonant terms, then this solvability condition is $F+iG = 0$.
We use \eqref{mrspertsol} and \eqref{mrspertsol1} in the right-hand side of \eqref{u2v2eq}, with $\xi$ or $\eta$ replaced by $z$, compute the coefficients of $e^{-iT}$, and impose the solvability condition. After some algebra, we find that $a(z,\tau)$ satisfies the degenerate Schr\"{o}dinger equation \eqref{asyeq}.

\subsection{Comparison with the MRS and linearized equations}\label{MRSlin}
Equation \eqref{asyeq} consists of the asymptotic equation \eqref{dqnls} for the MRS equation \eqref{gmrs} with a lower-order dispersive term.

When the MRS equation is derived from gas dynamics, with $u$ and $v$ the Eulerian velocity perturbations
in the right and left moving sound waves, there are additional
factors of $(\gamma+1)/2$ on the nonlinear terms with opposite signs \cite{mr}. This gives a factor on the quasilinear Schr\"odinger term $(|a|^2a_x)_x$
in \eqref{dqnls} of
\[
-\frac{2}{3} \left(\frac{\gamma+1}{2}\right)^2 = -\frac{(\gamma+1)^2}{6},
\]
which agrees with the one in  \eqref{asyeq} found by expanding the gas dynamics equations directly.

The linearization of the Lagrangian wave equation \eqref{gasdynamicsPDE} with sound speed \eqref{c0exp} is
\begin{equation}
w_{tt} - \left[1+\epsilon S(2x)\right] w_{xx} = 0.
\label{linwave}
\end{equation}
We suppose that $S(x)$ is a general real-valued, $2\pi$-periodic function with zero mean and Fourier expansion
\[
S(x) = \sum_{k\in \Z_*} S_k e^{ikx},
\]
and look for spatially periodic, resonantly reflected solutions of \eqref{linwave} of the form
\[
w =\left[A e^{ik x} + B e^{-ik x}\right] e^{-i\omega(k;\epsilon) t} + O(\epsilon),
\]
with wavenumber $k\in \Z_*$.

A perturbation expansion similar to the one above for the nonlinear equation, whose details we omit, shows that (up to an overall sign)
\begin{align*}
\omega(k;\epsilon) &= k +\epsilon\omega_1(k) + \epsilon^2 \omega_2(k) + O(\epsilon^3),
\\
\omega_1 &= \pm \frac{1}{2}k |S_k|,
\\
\omega_2 &= - \frac{\omega_1^2}{2k} - \frac{k}{16} \sum_{n\in \Z\setminus\{0,k\}} \frac{(2n-k)^2}{n(n-k)}\left|S_{n} \pm e^{i\sigma}S_{n-k}\right|^2,
\end{align*}
where $S_k = |S_k|e^{i\sigma}$ and $A = \pm e^{i\sigma}B$.

For a sawtooth wave with $S_k = 2i/k$ and $e^{i\sigma} = i\sgn k$, we find that
\begin{align*}
\omega_1 &= \pm 1,
\qquad
\omega_2 = - \frac{1}{2k}- \frac{1}{2}\pi^2 k.
\end{align*}
The first order perturbation $\omega_1$ is independent of the wavenumber $k$, corresponding to the nondispersive,
constant frequency oscillations in the resonant reflection of sound waves off a sawtooth entropy wave. (In contrast, other entropy profiles lead to dispersive oscillations.) The second order correction $\omega_2$ is the dispersion relation of
\[
a_{\tau} + \frac{1}{2}\partial_z^{-1}a - \frac{1}{2} \pi^2 a_z = 0,
\]
which corresponds to the  linear terms in \eqref{asyeq}.

We remark that the inclusion of higher order corrections in the sound speed perturbation $S$, say
\[
S(x, \epsilon) = S(x) + \epsilon R(x) + O(\epsilon^2),
\]
gives
\begin{align*}
\omega_1 &= \pm \frac{1}{2}k \left|S_{k} + \epsilon R_k + O(\epsilon^2)\right|
\\
&= \pm \frac{1}{2}k \left|S_k\right|\left(1 +  \epsilon q_k\right) + O(\epsilon^2),
\\
q_k &= \frac{1}{2}\left(\frac{R_k}{S_k}
+ \frac{R^*_k}{S^*_k}\right),
\end{align*}
which would lead to additional dispersive terms in \eqref{asyeq} if $q_k\ne 0$. Since the sawtooth wave $S(x)$ is odd and $S_k$ is imaginary, no such terms appear when $R(x)$ is even and $R_k$ is real. In particular, this is the case for the sound speed perturbations corresponding to a pure sawtooth entropy wave of the form
\[
s(x) =  s_0 + \epsilon S(2x),
\]
instead of \eqref{entexp}, so we would get an asymptotic equation of the same form as \eqref{asyeq} in that case also.

\section{Numerical solutions}\label{sec:num}

In this section, we show two sets of numerical solutions, one for front propagation in the MRS and DQS equations, the other for periodic solutions of the asymptotic gas dynamics equation.

\subsection{Front propagation for the MRS equation}\label{sec:num_mrs}

\begin{figure}
\centerline{
\includegraphics[width=0.5\textwidth]{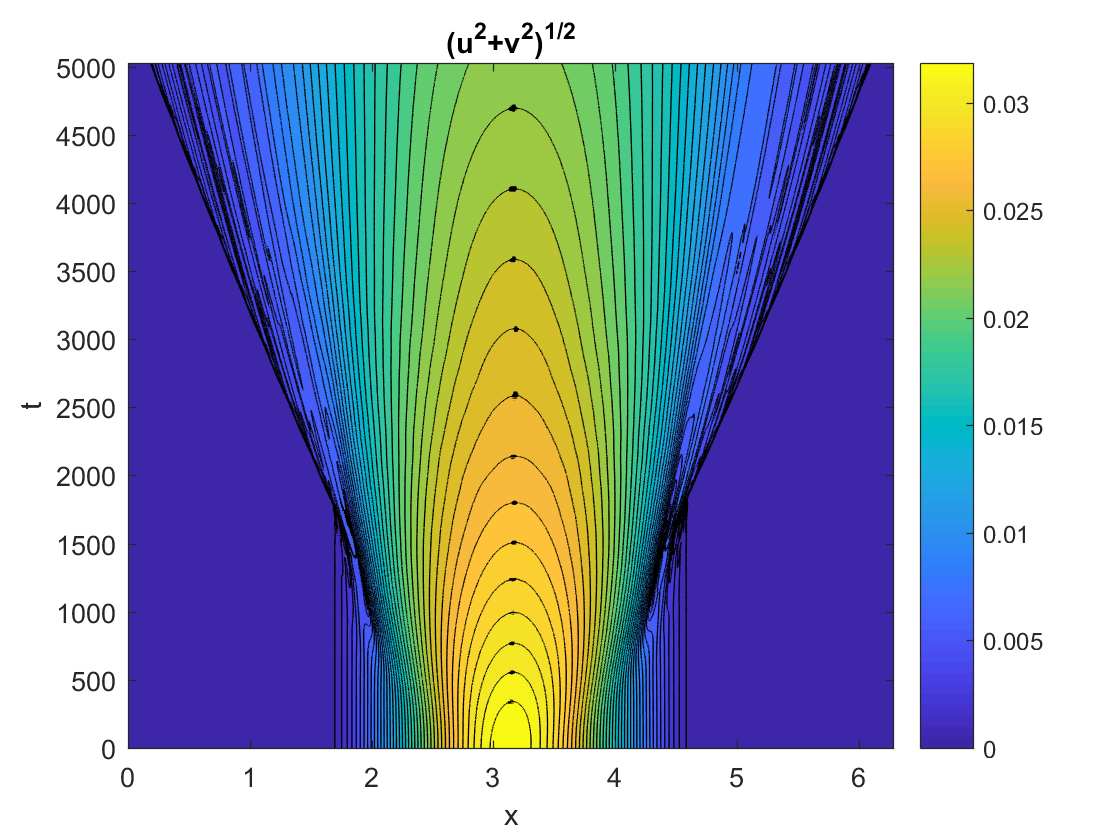}
\includegraphics[width=0.5\textwidth]{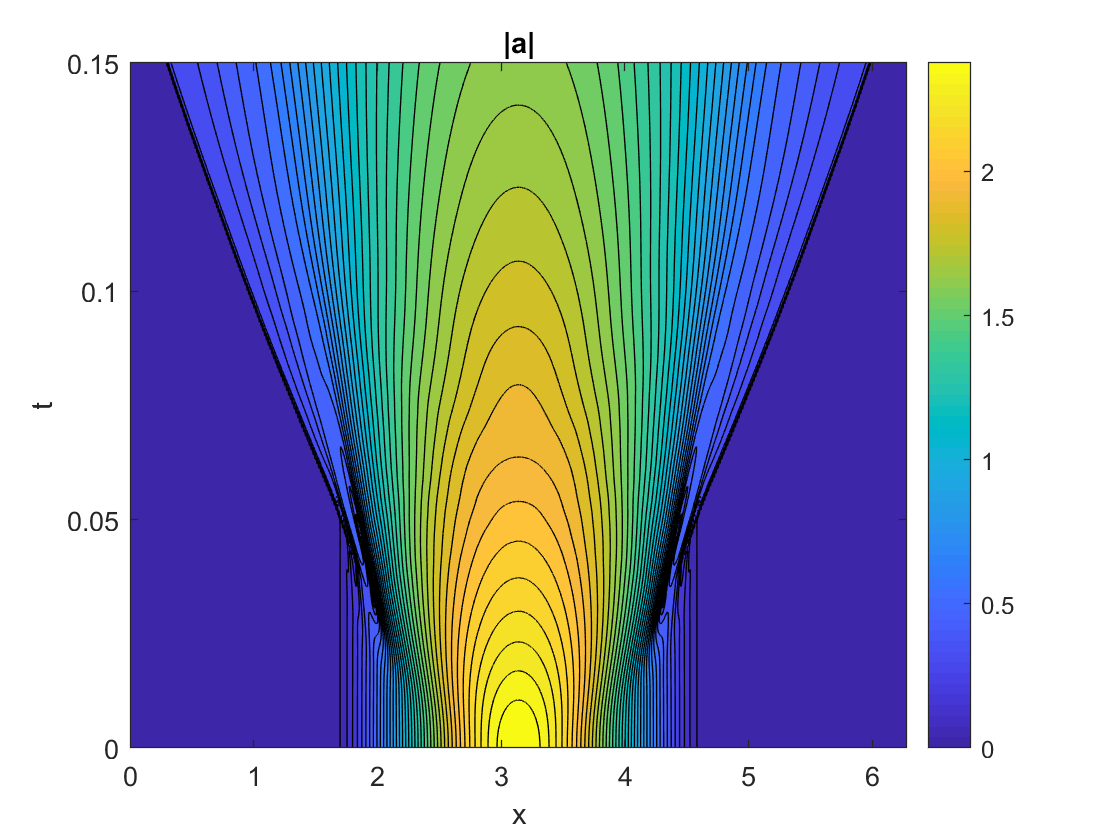}
}
\caption{Left: Contour plot of $(u^2+v^2)^{1/2}$ for the solution of the MRS equation \eqref{MRS} with initial data \eqref{mrsic}, where  $\epsilon = 0.0055$, and $0\le t \le 800\cdot 2\pi$. Right: Contour plot of the amplitude $|a|$, rescaled to the MRS solution, for the solution of the DQS equation \eqref{dqnls} on the corresponding slow time interval $0\le \tau \le 0.15$.}
\label{fig:contour}
\end{figure}

We consider the MRS equation \eqref{MRS} with the compactly supported initial data
\begin{align}
\begin{split}
u(x,0) &=
\begin{cases}
        \epsilon \left[\pi^2/4 - (x-\pi)^2\right]^2 &\text{if $|x-\pi| <\pi/2$},
        \\
        0    & \text{if $\pi/2 < |x-\pi| < \pi$},
\end{cases}
\\
 v(x,0) &= 0,
 \end{split}
 \label{mrsic}
\end{align}
where $\epsilon^2\cdot 800\cdot 2\pi = 0.15$, or $\epsilon \approx 0.005463$.
The left-hand side of Figure~\ref{fig:contour} shows a contour plot of $(u^2+v^2)^{1/2}$ for the solution of the MRS equation
The right-hand side of Figure~\ref{fig:contour} shows a contour plot of the corresponding solution of the DQS equation \eqref{dqnls} rescaled to the MRS variables.

The numerical solutions were computed by the method of lines, using an explicit forth order Runge-Kutta method in time and a forth order WENO flux in space for the MRS equation \cite{weno} with $2^{17}$ spatial grid points, and an implicit second order method in time and a pseudo-spectral method in space for the DQS equation with $2^{15}$ Fourier modes.
The pseudo-spectral method includes second-order spectral viscosity, which is required in order to get a convergent numerical solution
for the DQS equation with this compactly supported initial data.

\begin{figure}
\centerline{
\includegraphics[height=0.35\textwidth]{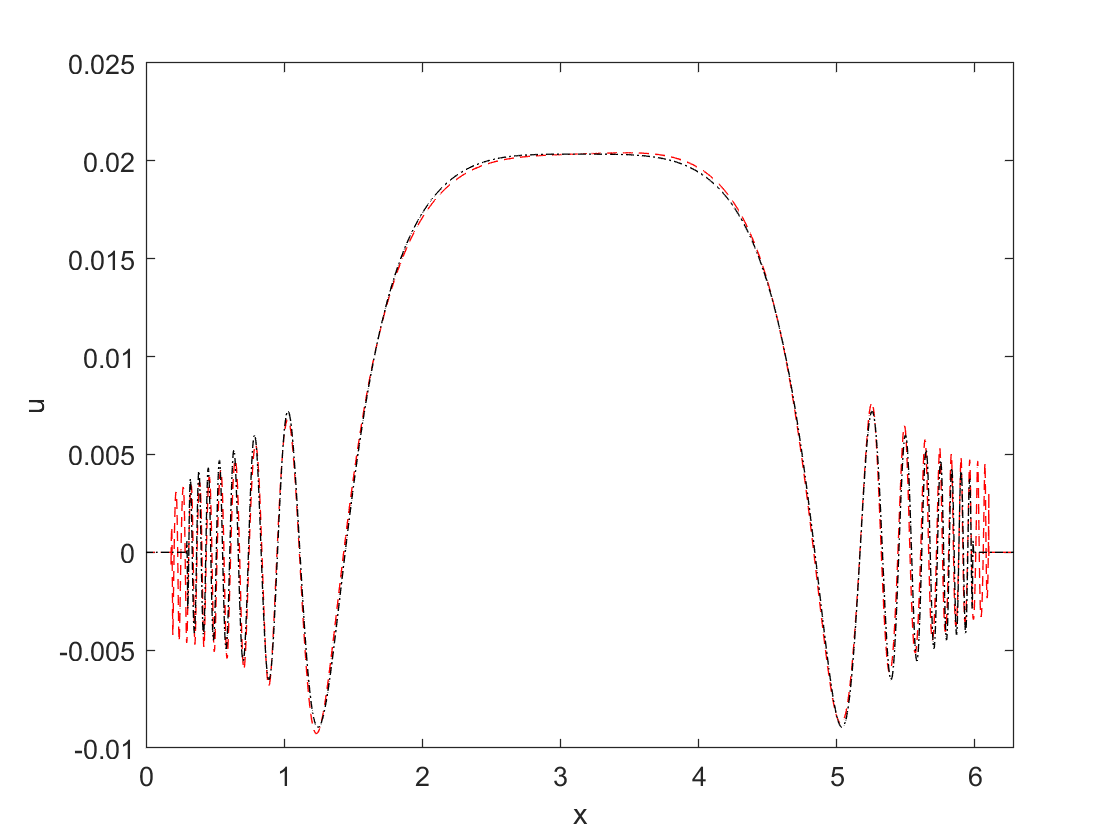}\quad
\includegraphics[height=0.35\textwidth]{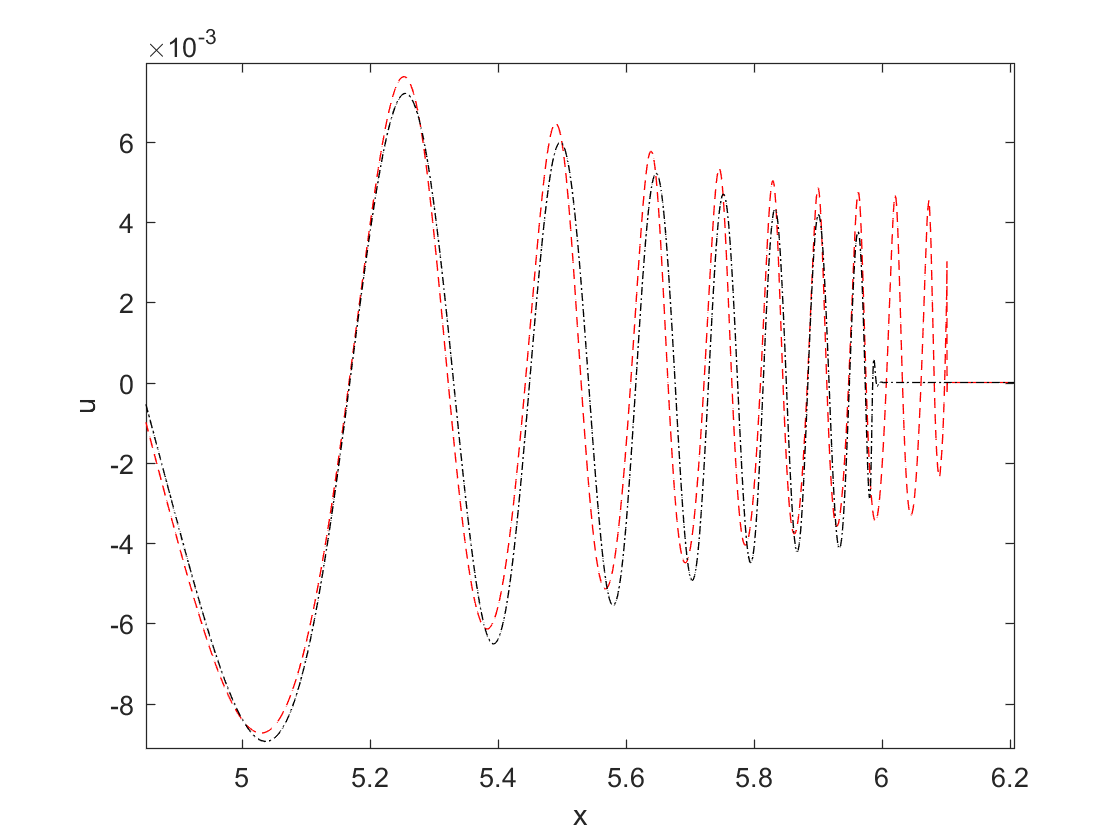}
}
\caption{Left: Solution on $0 <x < 2\pi$. Right: Detail near the front. The red line is the solution of the MRS equation \eqref{MRS}
 for $u(x,t)$ with initial data \eqref{mrsic} at $t=800\cdot 2\pi$. The black line is the corresponding
 asymptotic solution for $u(x,t)$ obtained from the DQS equation \eqref{dqnls} and
rescaled to the MRS variables.}
\label{fig:detail}
\end{figure}

The numerical solution of the MRS equation shows that the fronts at $x=\pi/2, 3\pi/2$ are initially stationary. A weak shock forms inside the left-hand side of the pulse at $t\approx 1050$, and a second shock forms closer to the front at $t\approx 1240$. This shock hits the left front at $t\approx 1760$, after which the front expands. A similar, but slightly later, sequence of events occurs at the right front. The MRS equation is not invariant under spatial reflections, rather it is invariant under
\[
x-\pi\mapsto -(x-\pi),\qquad  u\mapsto -v,\qquad v\mapsto -u,
\]
so complete reflectional symmetry is not to be expected. The mechanism of front expansion for the MRS equation is the propagation of shocks in $u$ or $v$ into the zero state.

The right-hand side of Figure~\ref{fig:contour} shows the corresponding numerical solution of the DQS equation. The solution agrees well with the MRS solution, although, as shown in Figure~\ref{fig:detail}, the DQS solution under-estimates the speed of the fronts. This discrepancy may be the result of greater numerical dissipation in the spectral viscosity scheme than in the WENO scheme. The DQS solution is symmetric in $x$, since averaging over the $(u,v)$-oscillations in the MRS solutions eliminates their reflectional asymmetry.

The structure of the DQS solution and the mechanism by which a DQS front propagates into the zero solution
is less clear than for the MRS equation.  In particular, compactly supported traveling waves or self-similar solutions
of the DQS equation \eqref{dqnls} that are dispersive analogs of the self-similar Barenblatt solutions for the porous medium equation are not weak solutions of \eqref{dqnls} at the front \cite{hs, smothers} and do not agree with the numerical simulations. These questions
will be analyzed further in \cite{hs}.

\subsection{Two-harmonic initial data for the asymptotic gas dynamics equation}\label{sec:num_gd}

First, we normalize the asymptotic gas dynamics equation \eqref{asyeq}. The dimensionless parameter $\epsilon$ in \eqref{c0exp} measures the strength of the perturbation in the sound speed
relative to the mean sound speed $\bar{c}$. A dimensionless parameter that measures the strength of the sound waves is $M = U/\bar{c}$ where $U$ is a typical size of velocity perturbations in the sound waves. In the asymptotic expansion, we assume that $M = O(\epsilon^{3/2})$.

We make the change of variables
\[
a = \frac{\sqrt{6} M}{(\gamma+1) \epsilon^{3/2}} \tilde{a},
\quad
\tilde{z} = z + \frac{1}{2}\pi^2 \tau,\quad \tilde{\tau} = \frac{M^2}{\epsilon^3} \tau,
\]
in \eqref{asyeq} and, after dropping the tildes, get
\begin{equation}
i\left(a_\tau + \mu \partial_z^{-1} a\right) = \left(|a|^2 a_z\right)_z,
\label{normasyeq}
\end{equation}
where $\mu = \epsilon^3/2M^2$ measures the relative strength of the dispersion and nonlinearity.
Up to an unimportant coefficient,
the dispersionless limit of this equation with $\mu=0$ is \eqref{dqnls}.

The solution of \eqref{normasyeq} with initial data $a(z,0) = a_0 e^{inz}$ is
\[
a(z,\tau) = a_0 e^{inz - i\Omega\tau},\qquad \Omega = -\frac{\mu}{n} - n^2|a_0|^2,
\]
where the only effect of nonlinearity is a frequency shift. Thus, we need at least two initial harmonics for the cubic nonlinearity to generate
a resonant cascade to higher harmonics and nontrivial dynamics.

\begin{figure}
\centerline{
\includegraphics[width=0.5\textwidth]{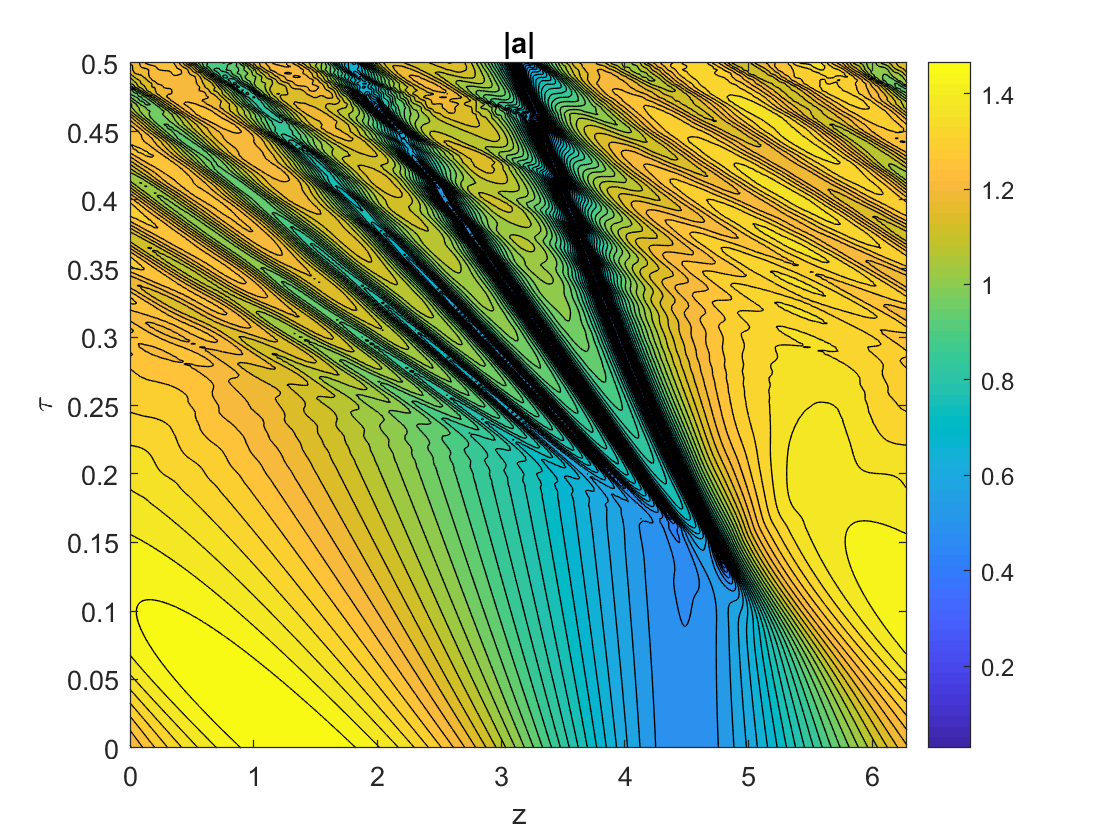}
\includegraphics[width=0.5\textwidth]{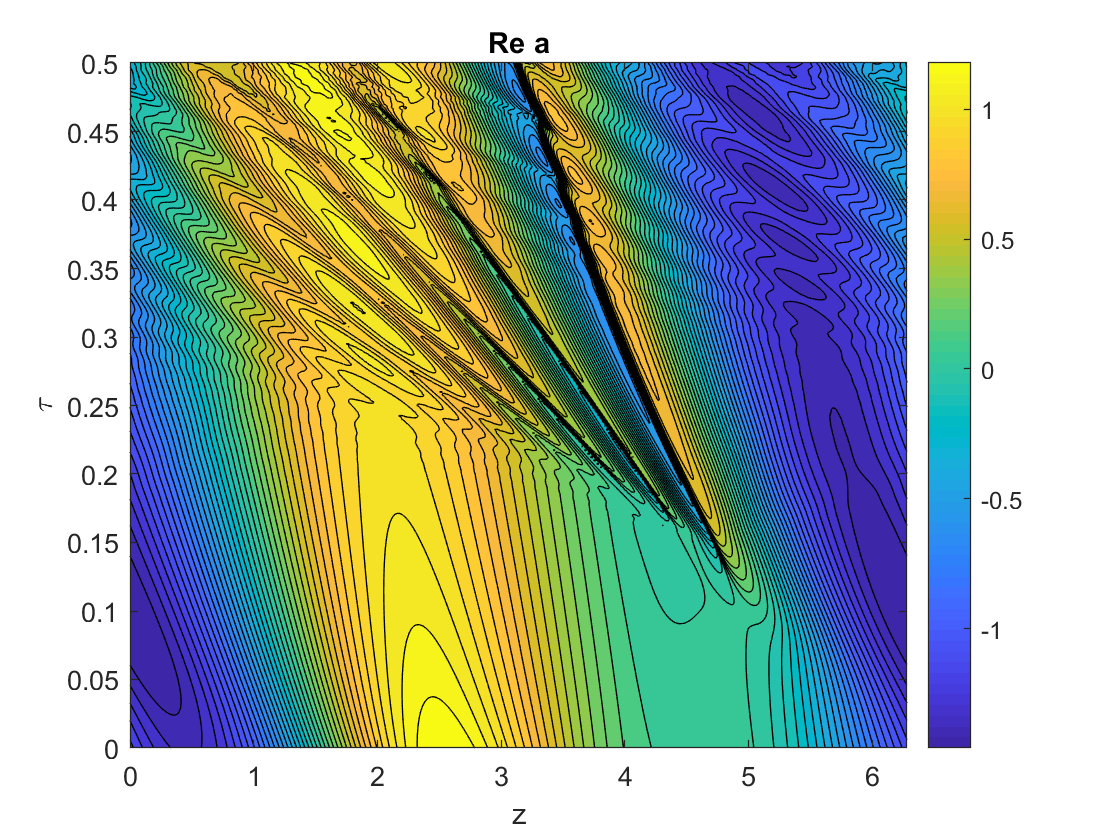}
}
\caption{A contour plot of the solution of \eqref{normasyeq} with $\mu=1$ and initial data \eqref{normic}. The solution is computed by a pseudo-spectral method using $2^{13}$ Fourier modes. Left: $|a|$. Right: $\re a$.}
\label{fig:2harmonics}
\end{figure}

\begin{figure}
\centerline{
\includegraphics[width=0.5\textwidth]{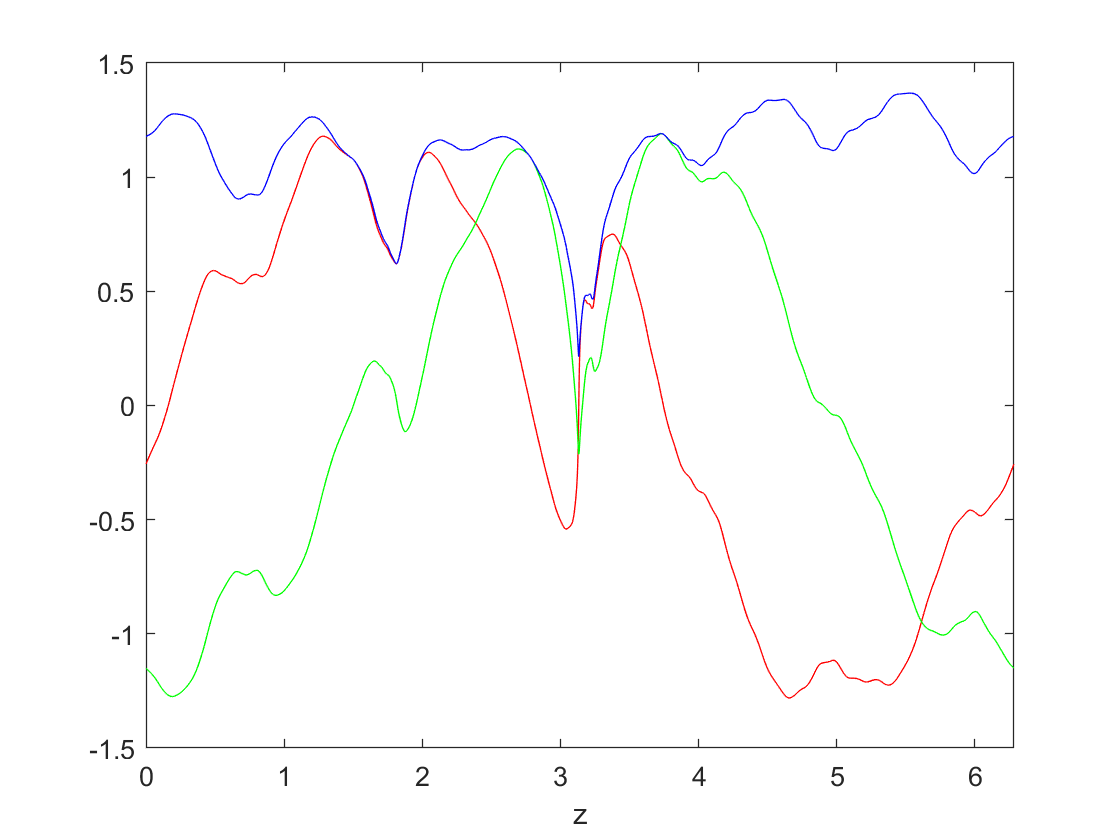}
\includegraphics[width=0.5\textwidth]{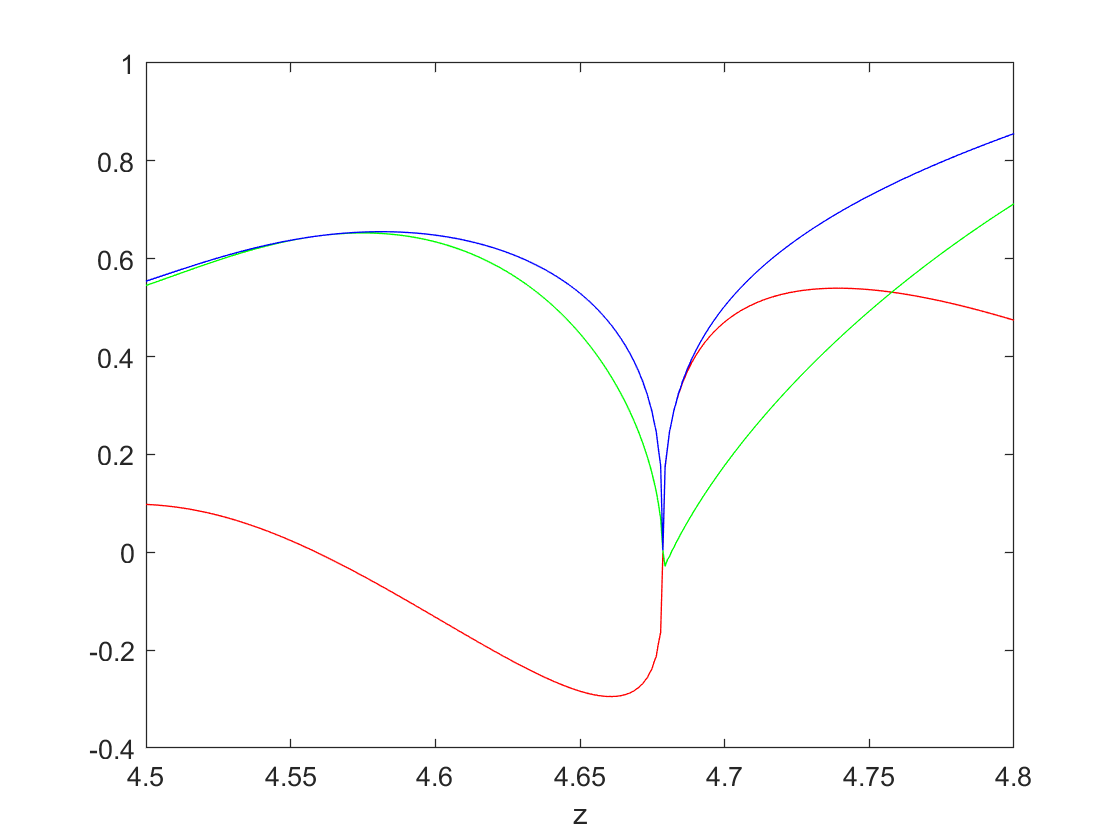}
}
\caption{Left: Graph of the solution of \eqref{normasyeq} with $\mu=1$ and initial data \eqref{normic} at $\tau =0.5$.
Right: Detail of the cusp in $|a|$ at $\tau=0.1581$. The minimum value of $|a|$ in the numerical solution is $|a| \approx 0.004$.
(Blue: $|a|$; Red: $\re a$; Green: $\im a$.)}
\label{fig:2harmonics_graph}
\end{figure}

In Figure~\ref{fig:2harmonics}, we show contour plots of $|a|$ and $\re a$ for the solution of \eqref{normasyeq} with $\mu=1$
and initial data
\begin{equation}
a(z,0) = - e^{iz} + \frac{1}{2} e^{2i(z + 2\pi^2)}.
\label{normic}
\end{equation}
In Figure~\ref{fig:2harmonics_graph}, we show a graph of the solution at time $\tau = 0.5$.
The solution is computed by the method of lines, using an implicit second order method in time and a
pseudo-spectral method in space. No spectral viscosity is required for the initial data \eqref{normic}.

The effects of dispersion on this solution are small. Numerical solutions for different values of $\mu$, including $\mu=0$, are qualitatively
similar to the one for $\mu=1$, and we do not show them here.

The numerical solution for $a$ becomes close to $0$ at $\tau\approx 0.15$, $z\approx 4.7$, when $|a|$ develops a cusp that propagates to the left; a graph of the solution near the cusp at $\tau = 0.1581$ is shown in Figure~\ref{fig:2harmonics_graph}. Further cusps form ahead of the original cusp, and these cusps later resolve into highly oscillatory waves in which $|a|$ is bounded well away from zero.

We remark that, although one cannot tell from the numerical solutions whether $|a|$ beomes exactly zero, the dispersionless version of \eqref{normasyeq}, with $\mu=0$, has traveling waves in which $a$ passes through zero \cite{hs}. After normalization, this traveling wave is given implicitly by
\[
a(z,\tau) = 1 - e^{-ic\Phi(z-c\tau)},\qquad  \Phi(\xi) - \frac{1}{c} \sin \left[c\Phi(\xi)\right] = \frac{1}{2}\xi.
\]
This solution has a cube-root dependence on $(z-c\tau)^{1/3}$ near $z=c\tau$, which is consistent with the qualitative  behavior of the
numerical solution, and suggests that the solution for $a$ may become H\"older continuous with exponent $1/3$ when $a$
vanishes at a point.

\end{document}